%% file: ABBR02.tex
\documentclass[a4paper,12pt,leqno,english]{article}
\usepackage[utf8]{inputenc}
\usepackage[T1]{fontenc}
\usepackage{babel}
\usepackage{amsmath}
\usepackage{amssymb}
\usepackage{enumitem}
\usepackage{graphicx}
\usepackage{hyperref}
\usepackage{tcolorbox}
\usepackage{pictex}
\usepackage{tikz} %

\input rs_macros.tex

\numberwithin{equation}{section}

\newtheorem{theorem+}           {Theorem}      [section]
\newtheorem{definition+}  [theorem+]  {Definition}
\newtheorem{lemma+}  [theorem+]  {Lemma}
\newtheorem{corollary+}  [theorem+]  {Corollary}
\newtheorem{proposition+}  [theorem+]  {Proposition}
\newtheorem{example+}  [theorem+]  {Example}
\newtheorem{problem+}  [theorem+]  {Problem}
\newtheorem{remark+}  [theorem+]  {Remark}

\newenvironment{theorem}{\begin{theorem+}\sl}{\end{theorem+}\rm}

\newenvironment{corollary}{\begin{corollary+}\sl}{\end{corollary+}\rm}

\newenvironment{example}{\begin{example+}\rm}{\end{example+}\rm}

\title{{\Large \bf 
Polynomials with exponents in  compact convex sets  and 
associated weighted extremal functions - \\
Characterization of polynomials by $L^2$-estimates }}
\author{{Benedikt Steinar Magnússon, Álfheiður Edda Sigurðardóttir,} \\
{Ragnar Sigurðsson and Bergur Snorrason}}
\date{{}} 

\begin{document}
\maketitle

\begin{abstract} \noindent
The main result of this paper is that an entire function $f$
that is in $L^2(\C^n,\psi)$ with respect to the 
weight  $\psi(z)=2mH_S(z)+\gamma\log(1+|z|^2)$ is  a
polynomial  
with exponents in $m\widehat S_\Gamma$.
Here $H_S$ is the
logarithmic supporting function of a compact convex  set $S\subset \R^n_+$
with $0\in S$, $\gamma\geq 0$ is small enough in terms of $m$, 
	and   $\widehat S_\Gamma$
is the hull of $S$  with respect to a certain cone 
$\Gamma$ depending on $S$, $m$ and $\gamma$. 
An  example showing that in general 
$\widehat S_\Gamma$ can not be replaced by $S$ is constructed.

\end{abstract}

\section{Introduction}
\label{sec:01}

\noindent 
Let  $S$  be a compact convex subset of $\R^n_+$
with $0\in S$.  We let  ${\mathcal P}^S_m(\C^n)$ denote the space 
of  all polynomials $p$ of the form 
$$
p(z)=\sum_{\alpha\in (mS)\cap \N^n} a_\alpha z^\alpha, \qquad z\in \C^n,
$$
and set 
${\mathcal P}^S(\C^n)=\cup_{m\in \N}{\mathcal P}^S_m(\C^n)$.
The \emph{$S$-Lelong class} $\L^S(\C^n)\subset \PSH(\C^n)$ 
is defined in the following way 
in terms of the \emph{supporting function} $\varphi_S$ of $S$, 
$$\varphi_S(\xi)=\sup_{s\in S}\scalar s\xi, \qquad \xi\in \R^n,
$$
and the map  $\Log\colon \C^{*n}\to \R^n$,
$$\Log\,  z=(\log|z_1|,\dots,\log|z_n|), \qquad z\in \C^{*n}.
$$
We define the \emph{logarithmic supporting function} 
$H_S\in \PSH(\C^n)$ of $S$ by 
$$H_S=\varphi_S\circ \Log$$
on $\C^{*n}$ and extend the definition to the coordinate 
hyperplanes  by the formula
$$
H_S(z)=\varlimsup_{\C^{*n}\ni w\to z}H_S(w), \qquad z\in \C^n\setminus
\C^{*n}.
$$
The $S$-Lelong class $\L^S(\C^n)$ 
is defined as the set of all $u\in \PSH(\C^n)$ satisfying a growth
estimate of the form $u\leq c_u+H_S$ for some constant $c_u$. 

Let $N^S_s=\{\xi\in \R^n\,;\, \scalar s\xi=\varphi_S(\xi)\}$ 
denote the \emph{ normal cone} of $S$ at the point $s$ and  define
for  every cone $\Gamma\subset (\R^n\setminus N^S_0)\cup\{0\}$  
with $\Gamma\neq \{0\}$  the   \emph{ $\Gamma$-hull of $S$} by
\begin{equation}
  \label{eq:1.1}
\widehat S_\Gamma
=\{x\in \R_+^ n\,;\, \scalar x\xi\leq \varphi_S(\xi),\  \xi\in \Gamma\}.   
\end{equation}
We say that $S$ is \emph{$\Gamma$-convex} if $S=\widehat S_\Gamma$
and that $S$ is a \emph{lower set} if for every $s\in S$
the cube $[0,s_1]\times\cdots\times [0,s_n]$ is contained 
in $S$. Every lower set is $\Gamma$-convex with repect to the cone
$\Gamma=\R^n_+$. (See \cite{MagSigSigSno:2023}, Section 4.)  
Our main result is:

\begin{theorem} \label{thm:1.1}
Let $S$ be a compact convex subset of $\R_+^n$ with  $0\in S$,
$m\in\N^*$, and $d_m=d(mS,\N^n\setminus mS)$ denote
the euclidean distance between 
the sets $mS$ and $\N^n\setminus mS$.
Let $f\in \OO(\C^n)$, assume that 
\begin{equation}
  \label{eq:1.2}
\int_{\C^n}|f|^2(1+|\zeta|^2)^{-\gamma}e^{-2mH_S} \,
d\lambda <+\infty
\end{equation}
for some $0\leq \gamma < d_m$. Let  $\Gamma$
be the cone consisting  of all $\xi$ such that 
the angle between the vectors ${\mathbf 1}=(1,\dots,1)$ and $\xi$ is 
$\leq \arccos(-(d_m-\gamma)/\sqrt n)$ and let
$\widehat S_\Gamma$ be  the hull of $S$ with respect to the cone
$\Gamma$
defined by (\ref{eq:1.1}).
Then $f\in {\mathcal P}_m^{\widehat S_\Gamma}(\C^n)$.
\end{theorem}

	\begin{figure}[h]
		\centering
		\begin{tikzpicture}
      scale = 1.1, every node/.style={transform shape}]
			\clip (-2.8, -2.8) rectangle (2.8, 2.8);

			\draw[fill = gray!20]
				(-10, 3) --
				(0, 0) --
				(10, -35) --
				(10, 10);


			\draw[->] (0, 0) -- (2, 2);
			\node at (1.4, 1.8) {\small $\mathbf{1}$};
			\node at (2.3, -0.8) {\small $\Gamma$};

			\draw (0, -10) -- (0, 10);
			\draw (-10, 0) -- (10, 0);

			\node at (1.2, 0.5) {\small $\theta_m$};
			\path[clip] (1, 1) -- (0, 0) -- (10, -35) -- cycle; 
			\draw (0, 0) circle (1);
		\end{tikzpicture}

{\small The cone $\Gamma$ has  opening angle 
$\theta_m = \arccos(-(d_m - \gamma)/\sqrt{n})$}
	\end{figure}

\noindent
Observe that the largest possible $d_m$ is $1$ and 
smallest possible $\gamma$ is $0$, which implies that
the largest possible opening angle of the cone $\Gamma$
is $\arccos(-1/\sqrt n)$. 
If $\gamma_0$ is the infimum of $\gamma$ such that 
\eqref{eq:1.2} holds, then $\Gamma(\gamma_0) = 
\cup_{\gamma>\gamma_0} \Gamma(\gamma)$.
Therefore $f$ is a polynomial with exponents in 
$m \widehat S_{\Gamma(\gamma_0)} = \cap_{\gamma>\gamma_0}
m\widehat S_{\Gamma(\gamma)}$.
We are interested in determining conditions on $S$
which ensure that $\widehat S_\Lambda =S$  for some cone 
$\Lambda\subseteq \Gamma$.

\begin{corollary} \label{cor:1.2}  The function $f$ in
Theorem~\ref{thm:1.1} is in  ${\mathcal P}_m^{S}(\C^n)$
in the cases:
\begin{enumerate}
\item [{\bf(i)}] $S$ is a lower set.
\item [{\bf(ii)}] $S=\widehat S_{\Lambda}$ for some cone 
$\Lambda$ contained in 
$\{\xi\in \R^n \,;\, \scalar{\mathbf 1}\xi\geq 0\}$.
\item [{\bf(iii)}]  $(mS)\cap \N^n=(m\widehat S_{\Gamma})\cap \N^n$. 
\end{enumerate}
\end{corollary}
\newpage

Bayraktar, Hussung, Levenberg, and Perera \cite{BayHusLevPer:2020},
stated in Proposition 4.3 that a function $f\in\OO(\C^n)$ satisfying
(\ref{eq:1.2}) is in ${\mathcal P}^S_m(\C^n)$ when $S$ is a 
polytope.  
This is true for example in the cases mentioned in 
Corollary \ref{cor:1.2}. However this is not correct in general
as seen by Example \ref{ex:4.1}.

For an exposition of the general background to the results 
in this paper see
\cite{MagSigSigSno:2023} and the references therein.

\subsection*{Acknowledgment}  
The results of this paper are a part of a research project, 
{\it Holomorphic Approximations and Pluripotential Theory},
with  project grant 
no.~207236-051 supported by the Icelandic Research Fund.
We would like to thank the Fund for its support and the 
Mathematics Division, Science Institute, University of Iceland,
for hosting the project.

\section{Taylor expansions of entire functions}
\label{sec:02}

Let $f\in \OO(\C^n)$ be an entire function with Taylor expansion
$f(z)=\sum_{\alpha\in \N^n}a_\alpha z^\alpha$ at the origin.
Let us now derive an estimate of the coefficients
$a_\alpha=f^{(\alpha)}(0)/\alpha!$.

\smallskip 
The Cauchy formula implies that  
for every polycircle $C_r=\{z\in \C^n\,;\, |z_j|=r_j\}$ with 
center $0$ and polyradius  $r\in \R_+^{*n}$ we have  
\begin{equation}
  \label{eq:2.1}
  a_\alpha =\dfrac{f^{(\alpha)}(0)}{\alpha !}
=\dfrac 1{(2\pi i)^n}\int_{C_r} \dfrac{f(\zeta)}{\zeta^\alpha}\cdot
\dfrac{d\zeta_1\cdots d\zeta_n}{\zeta_1\cdots\zeta_n}.
\end{equation}
We parametrize $C_r$ by $[-\pi,\pi]^n\ni \theta\mapsto 
(r_1e^{i\theta_1},\dots,r_ne^{i\theta_n})$ and get
\begin{equation}
  \label{eq:2.2}
  a_\alpha 
=\dfrac 1{(2\pi)^n}\int_{[-\pi,\pi]^n} 
\dfrac{f(r_1e^{i\theta_1},\dots,r_ne^{i\theta_n})}{r^\alpha
e^{i\scalar\alpha\theta}}\, d\theta_1\cdots d\theta_n.
\end{equation}
We take $\sigma,\tau \in \R^n$ with $\sigma_j<\tau_j$ for 
$j=1,\dots,n$ and define 
\begin{gather}
  \label{eq:2.3}
  A_{\sigma,\tau}=\{\zeta\in \C^n\,;\, e^{\sigma_j}\leq
  |\zeta_j|<e^{\tau_j}\}\subset \C^n,\\
L_{\sigma,\tau}=\prod_{j=1}^n\big([e^{\sigma_j},e^{\tau_j}]
\times[-\pi,\pi]\big)\subset
\R^{2n}=\big(\R^2\big)^n, \quad \text{ and }  \nonumber \\ 
K_{\sigma,\tau}=\prod_{j=1}^n [\sigma_j,\tau_j] \subset \R^n. \nonumber
\end{gather} 
We multiply both sides of (\ref{eq:2.2})
by $r_1\cdots r_n\, dr_1\cdots dr_n$, integrate with 
respect to the variables $r_j$ over 
$[e^{\sigma_j},e^{\tau_j}]$, observe that
$\int_{[e^{\sigma_j},e^{\tau_j}]}r_j\, dr_j=\tfrac
12(e^{2\tau_j}-e^{2\sigma_j})$, and get
\begin{equation}
  \label{eq:2.4}
  a_\alpha 
=\dfrac 1{v(A_{\sigma,\tau})}
\int_{L_{\sigma,\tau}} 
\dfrac{f(r_1e^{i\theta_1},\dots,r_ne^{i\theta_n})}{r^\alpha
e^{i\scalar\alpha\theta}}\, (r_1\, dr_1d\theta_1)\cdots (r_n\,
dr_nd\theta_n),
\end{equation}
where
$v(A_{\sigma,\tau})=\pi^n\prod_{j=1}^n(e^{2\tau_j}-e^{2\sigma_j})$
is the volume of the polyannulus.  By switching back to the 
original coordinates $\zeta_j=r_je^{i\theta_j}$ we get 
\begin{equation}
  \label{eq:2.5}
a_\alpha=\dfrac 1{v(A_{\sigma,\tau})} \int_{A_{\sigma,\tau}}
\dfrac{f(\zeta)}{\zeta^\alpha} \, d\lambda(\zeta).  
\end{equation}
Now we assume that 
$f\in L^2(\C^n,\psi)\cap\OO(\C^n)$ for some
measurable $\psi\colon \C^n\to \overline \R$, i.e., 
\begin{equation}
  \label{eq:2.6}
   \|f\|_\psi^2=\int_{\C^n}|f|^2e^{-\psi}\, d\lambda <+\infty,
\end{equation}
By (\ref{eq:2.5}) and the Cauchy-Schwarz inequality 
\begin{align}
  \label{eq:2.7}
  |a_\alpha|&\leq 
\dfrac 1{v(A_{\sigma,\tau})} \int_{A_{\sigma,\tau}}
|f(\zeta)|e^{-\psi(\zeta)/2} \cdot\dfrac{e^{\psi(\zeta)/2}}{|\zeta^\alpha|}
 \, d\lambda(\zeta)\\
&\leq \dfrac {\|f\|_\psi}{v(A_{\sigma,\tau})}
\bigg(\int_{A_{\sigma,\tau}} 
\dfrac{e^{\psi(\zeta)}}{|\zeta^\alpha|^2}
 \, d\lambda(\zeta)\bigg)^{1/2}_.  
\nonumber
\end{align}
If  $\psi$ is rotationally invariant in each
variable $\zeta_j$, i.e.,
$\psi(r_1e^{i\theta_1},\dots,r_ne^{i\theta_n})=
\psi(r_1,\dots,r_n)$,  then we  introduce
logarithmic coordinates $\xi_j=\log r_j$ and set 
$\chi(\xi)=\tfrac 12 \psi(e^{\xi_1},\dots,e^{\xi_n})$, and observe that
$r_j\, dr_j=e^{2\xi_j}\, d\xi_j$ and get 
\begin{equation}
  \label{eq:2.8}
\int_{A_{\sigma,\tau}} 
\dfrac{e^{\psi(\zeta)}}{|\zeta^\alpha|^2}
 \, d\lambda(\zeta)  
=
(2\pi)^n
\int_{K_{\sigma,\tau}} 
e^{2(\chi(\xi)-\scalar \alpha\xi+\scalar{\mathbf 1}\xi)} \, d\lambda(\xi).    
\end{equation}
We have 
$v(A_{\sigma,\tau})=\pi^n\prod_{j=1}^n(e^{2\tau_j}-e^{2\sigma_j})$, so
by combining  (\ref{eq:2.7}) and  (\ref{eq:2.8}) we get 
\begin{equation}
  \label{eq:2.9}
  |a_\alpha| \leq 
  \dfrac {\|f\|_\psi}{\prod_{j=1}^n(1-e^{-2(\tau_j-\sigma_j)})}
e^{-2\scalar{\mathbf 1}\tau}
\bigg(
\int_{K_{\sigma,\tau}} 
e^{2(\chi(\xi)-\scalar \alpha\xi+\scalar{\mathbf 1}\xi)} \, 
d\lambda(\xi)\bigg)^{1/2}_.    
\end{equation}
Since $\xi_j\leq \tau_j$ for every $\xi\in K_{\sigma,\tau}$ we arrive
at the estimate
\begin{equation}
  \label{eq:2.10}
  |a_\alpha| \leq 
  \dfrac {\|f\|_\psi}{\prod_{j=1}^n(1-e^{-2(\tau_j-\sigma_j)})}
e^{-\scalar{\mathbf 1}\tau}
\bigg(\int_{K_{\sigma,\tau}} 
e^{2(\chi(\xi)-\scalar \alpha\xi)} \, d\lambda(\xi)\bigg)^{1/2}_.    
\end{equation}

\section{Proof of Theorem \ref{thm:1.1}}
\label{sec:03}

\noindent
Let $f(z)=\sum_{\alpha\in \N^n} a_\alpha z^\alpha$ be
  the Taylor expansion of $f$ at the origin.  
We need to show
that $a_\alpha=0$ for every $\alpha\in \N^n\setminus m\widehat S_\Gamma$.
Since $\alpha\not\in m \widehat S_\Gamma$, there exists $\tau \in \Gamma$
such that $|\tau|=1$ and $\scalar \alpha\tau >m\varphi_S(\tau)$.
By rotating $\tau$ we may assume that $\tau$ is an interior point
of $\Gamma$ which gives $-\scalar{\mathbf 1}\tau <d_m-\gamma_0$.
We choose $\varepsilon >0$ such that
 $d_m-\gamma-\varepsilon>0$, and
$-\scalar{\mathbf 1}\tau <d_m-\gamma-\varepsilon$.

Recall that $\scalar\alpha\tau-m\varphi_S(\tau)$ is the euclidean
distance from $\alpha$ to the supporting hyperplane
$\{x\,;\, \scalar x\tau=m\varphi_S(\tau)\}$, so by assumption
$m\varphi_S(\tau)-\scalar\alpha\tau\leq -d_m$.   Hence
\begin{equation*}
-\scalar{\mathbf 1}\tau + m\varphi_S(\tau)-\scalar\alpha\tau
<-\gamma-\varepsilon.
\end{equation*}
We choose $\sigma\in \R^n\setminus \{0\}$ such that $\sigma_j<\tau_j$ for
every $j=1,\dots,n$  and 
\begin{equation*}
-\scalar{\mathbf 1}\tau + m\varphi_S(\xi) -\scalar\alpha\xi
<-(\gamma+\varepsilon)|\xi|, \qquad \xi\in K_{\sigma,\tau}.
\end{equation*}
By homogeneity we get 
\begin{equation}
  \label{eq:3.1}
-t\scalar{\mathbf 1}\tau + m\varphi_S(\xi)-\scalar\alpha\xi
<-(\gamma+\varepsilon)|\xi|, \qquad t>0, \ \xi\in tK_{\sigma,\tau}.
\end{equation}
Let $\xi_j=\log|\zeta_j|$ and observe that
$(1+|\zeta|^2)^\gamma \leq (n+1)^\gamma
\max\{1,\|\zeta\|_\infty^{2\gamma} \}$. From this inequality and
\eqref{eq:1.2} it follows that
$f\in L^2(\C^n,\psi)$,  where
\begin{equation}
  \label{eq:3.2}
\tfrac 12 \psi(\zeta)={\gamma}\log\|\zeta\|_\infty+mH_S(\zeta)  
=\gamma \|\xi\|_\infty+m\varphi_S(\xi).
\end{equation}
We set $\chi(\xi)=\tfrac 12 \psi(e^{\xi_1},\dots,e^{\xi_n})$.
Then the estimate (\ref{eq:3.1}) gives  
\begin{equation*}
-t\scalar{\mathbf 1}\tau + \chi(\xi)-\scalar\alpha\xi
<-\varepsilon|\xi|, \qquad t>0, \ \xi\in tK_{\sigma,\tau},
\end{equation*}
the estimate (\ref{eq:2.10}) with $tK_{\sigma,\tau}$ in the
role  of $K_{\sigma,\tau}$ gives 
\begin{align*}
  |a_\alpha| &\leq 
  \dfrac {\|f\|_\psi}{\prod_{j=1}^n(1-e^{-2(\tau_j-\sigma_j)t})}
e^{-t\scalar{\mathbf 1}\tau}
\bigg(\int_{tK_{\sigma,\tau}} 
e^{2(\chi(\xi)-\scalar\alpha\xi)} \, d\lambda(\xi)\bigg)^{1/2}\\
&\leq \dfrac {\|f\|_\psi}{\prod_{j=1}^n(1-e^{-2(\tau_j-\sigma_j)t})}
e^{-\epsilon|\sigma|t}t^{n/2} v(K_{\sigma,\tau})^{1/2} \to 0,
\qquad t\to +\infty,
\nonumber    
\end{align*}
and we conclude that $a_\alpha=0$.

\section{An example}
\label{sec:04}

Observe that (\ref{eq:1.2}) is equivalent to $\|f\|_\psi<+\infty$ 
in the notation presented in  \eqref{eq:2.6} for
$\psi\in \mathcal{PSH}(\C^n)$ given by 
$$\psi(z)= 2m H_S(z)+\gamma \log (1+|z|^2), \qquad z\in \C^n.$$
Note that $\|f\|_{2mH_S}\geq \|f\|_{\psi}$ for any $\gamma\geq 0$, so
if $\|f\|_{2mH_S}<+\infty$, then $\|f\|_{\psi}<+\infty$.

Now we are going to give an example of a compact convex neighbourhood
of 0 in $\R^n_+$ denoted by $S$, such that for every $m\geq 4$, there
exist polynomials $p$ such that
$\|p\|_{2mH_S}<+\infty$, but $p\notin \mathcal P^S_m(\C^n)$. This
shows that $\widehat S_\Gamma$ in Theorem \ref{thm:1.1} can not be
replaced by $S$.

\begin{example}\label{ex:4.1} 

Let $m\geq 4$ and $S\subseteq \R_+^2$ be the quadrilateral 
\begin{equation*}
S=\ch\{(0,0), (a,0), (b,1-b),(0,1)\}.  
\end{equation*}
where $0<a<1/m$ , $0<a<b<1$, $m(1-b)<1$, and $(b-a)/(1-b)>m-2-am$.

\smallskip
Then $(1,0),(2,0), \dots , (m-3,0)\notin mS$, but the following
calculations show that $\|p_k\|_{2mH_S}<+\infty$ for 
$$
p_k(z)=z^{(k,0)}=z_1^k, \qquad k=1,\dots, m-3.
$$

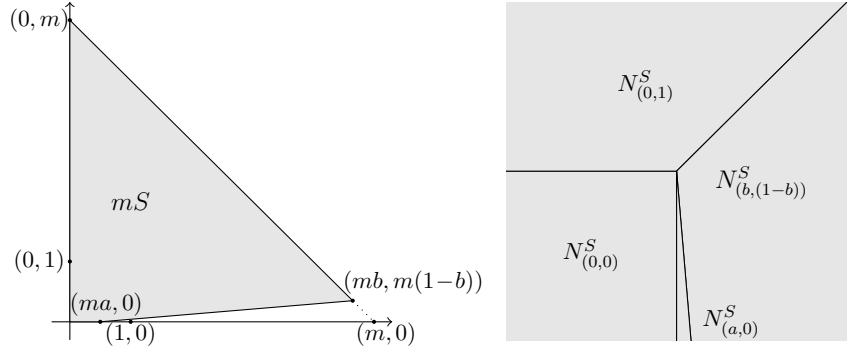
\begin{figure}[!h]
\def\mynda{0.1}
\def\myndb{0.93}
\def\myndm{5}
\def\myndscale{0.8} 
\centering
\begin{tikzpicture}[baseline={(2,2)},
scale = \myndscale, every node/.style={transform shape}]
	\clip (-1, -1) rectangle (7, 7);
	\draw[->] (-0.3, 0) -- (5.3, 0);
	\draw[->] (0, -0.3) -- (0, 5.3);

	\draw[fill = gray!20]
		(0, \myndm) --
		(\myndm*\myndb, \myndm - \myndm*\myndb) --
		(\myndm*\mynda, 0) --
		(0, 0) --
		(0, \myndm);
	\node at (1, 2) {$mS$};

	\draw[dotted] (\myndm*\myndb, \myndm - \myndm*\myndb) -- (\myndm, 0);

	\node[draw, fill, circle, inner sep = 0.5pt] at (0, \myndm) {};
	\node at (-0.5, \myndm) {\small $(0, m)$};

	\node[draw, fill, circle, inner sep = 0.5pt] at (\myndm*\mynda, 0) {};
	\node at (\myndm*\mynda +0.1, 0.3) {\small $(ma, 0)$};


	\node[draw, fill, circle, inner sep = 0.5pt] at (\myndm, 0) {};
	\node at (\myndm + 0.2, -0.2) {\small $(m, 0)$};

	\node[draw, fill, circle, inner sep = 0.5pt] at (0, 1) {};
	\node at (-0.5, 1) {\small $(0, 1)$};

	\node[draw, fill, circle, inner sep = 0.5pt] at (1, 0) {};
	\node at (1, -0.2) {\small $(1, 0)$};

	\node[draw, fill, circle, inner sep = 0.5pt] at (\myndm*\myndb, \myndm - \myndm*\myndb) {};
	\node at (\myndm*\myndb + 1, \myndm - \myndm*\myndb + 0.3) {\small $(mb, m(1\! -\! b))$};
\end{tikzpicture}
\begin{tikzpicture}[baseline={(0,0)},
scale = \myndscale, every node/.style={transform shape}]
	\clip (-2.8, -2.8) rectangle (2.8, 2.8);
	\draw[fill = gray!20]
		(0, 0) --
		(-10, 0) --
		(0, -10) --
		(0, 0);
	\node at (-1.4, -1.4) {\small $N^{S}_{(0,0)}$};

	\draw[fill = gray!20]
		(0, 0) --
		(-5, 0) --
		(-5, 5) --
		(5, 5) --
		(0, 0);
	\node at (-0.5, 1.4) {\small $N^{S}_{(0,1)}$};

	\draw[fill = gray!20]
		(0, 0) --
		(1 , {\myndm - \myndb*\myndm - 1*(\myndb - \mynda)/(1 - \myndb) - \myndm + \myndb*\myndm}) --
		(10, \myndm - \myndb*\myndm + 10 - \myndm + \myndb*\myndm) --
		(0, \myndm - \myndb*\myndm - \myndm + \myndb*\myndm);
	\node at (1.4, -0.2) {$N^{S}_{(b,(1-b))}$};

	\draw[fill = gray!20]
		(0, 0) --
		(1, {-1*(\myndb - \mynda)/(1 - \myndb)}) --
		(0, -4) --
		(0, 0);
	\node at (0.9, -2.5) {\small $N^{S}_{(a,0)}$};
\end{tikzpicture}
\end{figure}

\noindent
Since the map $(\R\times ]-\pi,\pi[)^2\to 
\C^2$, $(\xi_1,\theta_1,\xi_2,\theta_2)\mapsto 
(e^{\xi_1+i\theta_1},e^{\xi_2+i\theta_2})$ has the Jacobi determinant
$e^{2\xi_1+2\xi_2}$, we have
\begin{align*}
  \|p_k\|_{2mH_S}&=\int_{\C^2}|z_1|^{2k}e^{-2mH_S(z)}\, d\lambda(z)
=4\pi^2 \int_{\R^2}e^{2(k+1)\xi_1+2\xi_2-2m\varphi_S(\xi)}\,
                 d\xi_1d\xi_2.
\end{align*}
Observe that from \cite{MagSigSigSno:2023}, equations (3.5) and (3.6), we get
\begin{equation*}
 \varphi_S(\xi)=
	\max_{x\in \operatorname{ext}(S)} \langle x,\xi \rangle
	=
 \begin{cases}
 0, &\xi \in N^S_{(0,0)},\\   
a\xi_1, &\xi \in N^S_{(a,0)},\\
b\xi_1+(1-b)\xi_2,  &\xi \in N^S_{(b,(1-b))},\\
\xi_2, &\xi \in N^S_{(0,1)}.  
 \end{cases}
\end{equation*}
\noindent
We split the  integral over $\R^2$ into the sum of the
integrals over the normal cones at the extreme points of $S$, which we
calculate as
\begin{align*}
  \int_{N^S_{(0,0)}} e^{2(k+1)\xi_1+2\xi_2}\, d\xi
&=\int_{-\infty}^0 e^{2(k+1)\xi_1} \, d\xi_1 \,   \int_{-\infty}^0
e^{2\xi_2} \, d\xi_2=\dfrac 1{4(k+1)},\\
  \int_{N^S_{(a,0)}} e^{2(k+1)\xi_1+2\xi_2-2ma\xi_1}\, d\xi
&=\int_0^{\infty} e^{2(k+1-ma)\xi_1}\,
    \int_{-\infty}^{-\xi_1(b-a)/(1-b)} e^{2\xi_2}\, d\xi_2\, d\xi_1\\
& =\dfrac {1}{4((b-a)/(1-b)+ma-1-k)},
\end{align*}

\vspace{-0.3cm}

\begin{multline*}
 \int_{N^S_{(b,1-b)}} e^{2(k+1)\xi_1+2\xi_2-2m(b\xi_1+(1-b)\xi_2)}\,
   d\xi \\ 
=\int_0^{\infty} e^{2(k+1-mb)\xi_1}\,
\int_{-\xi_1(b-a)/(1-b)}^{\xi_1} e^{2(1-m(1-b))\xi_2}\, d\xi_2\,
                                  d\xi_1 \\
 =\dfrac {1}{4(1-m(1-b))}\bigg(
\dfrac 1{m-2-k}+\dfrac{1}{(b-a)/(1-b)+ma-1-k)}
\bigg),
\end{multline*}

\vspace{-0.2cm}

\begin{align*}
 \int_{N^S_{(0,1)}} e^{2(k+1)\xi_1+2\xi_2-2m\xi_2}\, d\xi\
&=\int_0^\infty e^{2(1-m)\xi_2}
\int_{-\infty}^{\xi_2}e^{2(k+1)\xi_1}\, d\xi_1 \, d\xi_2\\
&=\dfrac {1}{4(k+1)(m-2-k)}.
\end{align*}

\noindent
This shows that $\|p_k\|_{2mH_S}<+\infty$, and we have found polynomials
satisfying (\ref{eq:1.2}) which are not in $\mathcal P^S_m(\C^n)$.
\end{example}

{\small 
\bibliographystyle{siam}
\bibliography{rs_bibref}

\smallskip\noindent
Science Institute, 
University of Iceland, 
IS-107 Reykjav\'ik, 
ICELAND 

\smallskip\noindent
bsm@hi.is, alfheidur@hi.is, ragnar@hi.is, bergur@hi.is.
}

\end{document}

%% file: rs_macros.tex
\newcommand{\C}{{\mathbb  C}}

\renewcommand{\L}{{\mathcal L}}

\newcommand{\N}{{\mathbb  N}}
\newcommand{\OO}{{\mathcal O}}

\newcommand{\R}{{\mathbb  R}}

\newcommand{\scalar}[2]{{\langle#1,#2\rangle}}

\newcommand{\ch}{{\operatorname{ch}}}

\newcommand{\Log}{{\operatorname{Log}}}

\newcommand{\PSH}{{\operatorname{{\mathcal{PSH}}}}}

\def\vfigura#1#2{
\setbox0\vbox{{
\input #1
}}
\setbox1\vbox{\hbox{\box0}\hbox{{\obeylines #2}}}
\dimen0 = -\ht1
\advance\dimen0 by-\dp1
\dimen1 = \wd1
\dimen2 = -\dimen0
\divide\dimen2 by\baselineskip
\count100 = 1
\advance\count100 by\dimen2
\advance\count100 by1
\box1
\hangindent\dimen1
\hangafter=-\count100
\vskip\dimen0
}